\documentclass[12pt]{article}
\usepackage{amssymb}
\newcommand{\be}{\begin{equation}}
\newcommand{\ee}{\end{equation}}
\newcommand{\bea}{\begin{eqnarray}}
\newcommand{\eea}{\end{eqnarray}}
\newcommand{\bez}{\begin{eqnarray*}}
\newcommand{\eez}{\end{eqnarray*}}

\renewcommand{\d}{{\rm d}}
\newcommand{\A}{{\cal A}}
\newcommand{\V}{{\cal V}}

\begin{document}

\title{Differential Calculi on Quantum Spaces \\ determined by Automorphisms}

\author{Aristophanes Dimakis \\
 Department of Financial and Management Engineering \\
 University of the Aegean, 31 Fostini Str., GR-82100 Chios, Greece \\
 E-mail: dimakis@aegean.gr
 \and
 Folkert M\"uller-Hoissen \\
 Max-Planck-Institut f\"ur Str\"omungsforschung \\
 Bunsenstrasse 10, D-37073 G\"ottingen, Germany \\
 E-mail: Folkert.MuellerHoissen@mpi-sf.mpg.de }

\date{ }

\maketitle

\renewcommand{\theequation} {\arabic{section}.\arabic{equation}}

\begin{abstract}
If the bimodule of 1-forms of a differential calculus over an associative
algebra is the direct sum of 1-dimensional bimodules, a relation
with automorphisms of the algebra shows up. This happens for some familiar
quantum space calculi.
\end{abstract}

\section{Introduction}
\setcounter{equation}{0}
A {\em differential calculus} over an {\em associative} algebra $\A$
is given by a graded associative algebra $\Omega(\A) = \bigoplus_{r \geq 0} \Omega^r$,
where $\Omega^0 = \A$, together with a linear map
$\d : \; \Omega^r \rightarrow \Omega^{r+1}$ satisfying $\d^2 = 0$
and the graded Leibniz rule. Choosing generators $\theta^s$, $s \in S$, of the
$\A$-bimodule $\Omega^1$ of 1-forms, the commutation relations with elements
$a \in \A$ can be expressed as
\be
   \theta^s \, a = \sum_{s' \in S} \Phi(a)^s_{s'} \; \theta^{s'} \qquad \forall a \in \A
\ee
where $\Phi$ is a homomorphism from $\A$ into the $|S|\times |S|$-matrices with entries
in $\A$. Suppose we can choose $\theta^s$ such that $\Phi$ becomes diagonal:
\be
   \theta^s \, a = \phi_s(a) \; \theta^s  \qquad \forall s \in S  \; .  \label{theta-phi}
\ee
Then the maps $\phi_s$, $s \in S$, have to be \emph{automorphisms} of $\A$. In this case,
$\Omega^1$ is the direct sum of 1-dimensional bimodules.
If the special form (\ref{theta-phi}) of the commutation relations can be achieved,
calculations in the differential algebra are drastically simplified.
Trivial examples are ordinary differential forms on a parallelizable manifold and
K\"ahler differentials on freely generated commutative algebras. In these cases,
we have $\phi_s = {\rm id}$, which means that the basis 1-forms commute with all
elements of $\A$. Surprisingly, this can also be achieved for certain differential
calculi on \emph{non}commutative algebras \cite{DMH93sgqm,Dima+Mado96}.
Simple examples with non-trivial automorphisms are given by $\d x \; x = (x + \ell) \, \d x$
\cite{DMH92qml,DMHS93a,DMHS93b} and $\d x \; x = q \, x \, \d x$ \cite{DMH92qml,Mani92,Karo01}
on the commutative algebra freely generated by $x$, setting $\theta^1 = \d x$.
A large class of examples arises from differential calculi on Cayley graphs (`group lattices')
of a finite group $G$, where $S \subset G \setminus \{e\}$ and $\phi_s = R_s^\ast$,
the pull-back with the right action of $s \in S$ on $G$ \cite{DMH03grpl1,DMH03grpl2}.
There are more examples \cite{DMH04auto}, some will be discussed in this work.
Section~\ref{sec:auto->diff} recalls a recipe \cite{DMH04auto} to construct differential
calculi with the property (\ref{theta-phi}). The two subsequent sections present some
`quantum space' examples.
Section~\ref{sec:diffgeo} shows how pseudo-Riemannian structures can be introduced on algebras
with such a differential calculus.

\section{From automorphisms to differential calculi}
\label{sec:auto->diff}
\setcounter{equation}{0}
Let $\A$ be an associative algebra and $\{ \phi_s \, | \, s \in S \}$ a finite set of
{\em automorphisms}. Naturally associated with the latter are the \emph{discrete derivatives}
\be
     e_s = \phi_s - \mbox{id} \; .
\ee
A generalization is given by {\em twisted inner derivations}:
\be
    e_s = a_s \, \phi_s - R_{a_s}   \label{twist-inner-deriv}
\ee
with special elements $a_s \in \A$ and the right action $R_a$ by $a$.
As a consequence,
\be
  e_s(a \, a') = e_s(a) \, \phi_s(a') + a \, e_s(a')
  \qquad  \forall \, a, a' \in \A
\ee
so that $e_s$ is a \emph{twisted derivation}.
Other twisted derivations are obtained from twisted inner derivations
as limits
\be
    e_s = \lim_{t \to 0} \frac{1}{t} ( a_s \, \phi_s^t - R_{a_s} )
\ee
where $\phi_s^t$ is a one-parameter group of automorphisms with $\phi_s^0 = \mbox{id}$.

The $\A$-bimodule $\mathbb{X}$ of `vector fields' is generated
by $\{ e_s \}$ such that
\be
    (e_s \cdot a) \, a' := e_s(a') \, a   \qquad
    a \cdot e_s := e_s \cdot \phi_s(a) \; .
\ee
In the following, we require that $\{ e_s \, | \, s \in S \}$ is a left and right
$\A$-module basis of $\mathbb{X}$. A dual $\A$-bimodule $\hat{\Omega}^1$ of `1-forms'
is then determined via
\be
       \langle a \, \omega \, , \, X \rangle
  := a \, \langle \omega \, , X \rangle  \, , \quad
      \langle \omega \, , \, X \cdot a \rangle
  := \langle \omega \, , X \rangle \, a  \, , \quad
      \langle \omega \, a \, , \, X \rangle
  := \langle \omega \, , \, a \cdot X \rangle
\ee
for all $X \in \mathbb{X}, \, \omega \in \hat{\Omega}^1$.
The basis $\{ \theta^s \, | \,  s \in S \}$ dual to $\{ e_s \}$, i.e.
$\langle \theta^s \, , \, e_{s'} \rangle = \delta^s_{s'}$,
then satisfies (\ref{theta-phi}).
Now $\langle \d a \, , \, X \rangle := X a$ determines a linear map
$\d : \A \rightarrow \hat{\Omega}^1$ such that
\be
    \d a = \sum_{s \in S} e_s(a) \; \theta^s  \; .  \label{da-e_s}
\ee
The map $\d$ satisfies the Leibniz rule $\d (a \, a') = (\d a) \, a' + a \, (\d a')$
if $e_s$, $s \in S$, are twisted derivations. Furthermore, $(\hat{\Omega}^1, \d)$ is {\em inner}:
\be
   \d a = [ \vartheta \, , \, a ] \quad \mbox{with} \quad
   \vartheta := \sum a_s \, \theta^s \; .   \label{vartheta}
\ee
\noindent
{\em Remark.} In general, $\hat{\Omega}^1$ is larger than the $\A$-bimodule
$\Omega^1$ generated by $\d \A$, since $\theta^s$ need not be of the form
$\sum_i a_i \, \d a_i'$. If there are `coordinates' $x^s \in \A$
s.t. $e_s(x^{s'})$ are the components of an invertible matrix, then
$\theta^s \in \Omega^1$ and the two spaces coincide.
In most cases, however, $\A$ can be appended by {\em additional relations}
(typically: {\em inverses} of some elements), such that
$\hat{\Omega}^1(\hat{\A}) = \Omega^1(\hat{\A})$
for the new algebra $\hat{\A}$.
\hfill $\blacksquare$
\vskip.1cm

The first order differential calculus extends to higher orders
such that
\be
  \d \omega = [ \vartheta \, , \, \omega ]_{\mathrm{grad}} - \Xi(\omega)
\ee
involving the graded commutator. The map $\Xi$ has to satisfy
\be
 \Xi(a \, \omega \, a') = a \, \Xi(\omega) \, a' \, , \quad
 \Xi( \omega \, \omega') = \Xi(\omega) \, \omega'
    + (-1)^{\mathrm{grad}(\omega)} \, \omega \, \Xi(\omega') \; .
\ee

\section{WZ calculus on the quantum plane}
\label{sec:qp}
\setcounter{equation}{0}
The $q$-plane (see \cite{Klim+Schm97}, for example) is the algebra generated by two
elements $x,y$ subject to $x \, y = q \, y \, x$ with $q \in \mathbb{C} \setminus \{ 0,\pm 1 \}$.
Let us choose the special automorphisms
\be
    \phi_1(x) = r^{-1} \, x \, , \quad  \phi_1(y) = r^{-1} \, y \, , \quad
    \phi_2(x) = x \, , \quad \phi_2(y) = r^{-1} \, y
\ee
with $r \not\in \{ 0,1 \}$.
Applying the construction of the preceding section, starting with the
\emph{discrete derivatives} $e_s := \phi_s - \mbox{id}$, leads to
\be
   \theta^1 = {1 \over 1-r} \, \d x \; x^{-1}  \, , \qquad
   \theta^2 = {1 \over 1-r} \, ( \d y \; y^{-1} - \d x \; x^{-1} )
\ee
where we need inverses of $x$ and $y$, thus passing over to the quantum torus.
 From (\ref{theta-phi}) we recover the Wess-Zumino calculus \cite{Wess+Zumi90}
on the $q$-plane:
\bea
   \begin{array}{ll}
   x \; \d x = p \, q \; \d x \; x \qquad &
   y \; \d x = p \; \d x \; y  \\
   y \; \d y = p \, q \; \d y \; y  &
   x \; \d y = q \; \d y \; x + ( p \, q - 1) \; \d x \; y
   \end{array}
\eea
where $p := r/q$.
Furthermore, $\theta^1 \, \theta^2 = - \theta^2 \, \theta^1$, $(\theta^1)^2 = 0 = (\theta^2)^2$,
$\d \omega = [ \vartheta \, , \, \omega ]_{\mathrm{grad}}$ (so that $\Xi = 0$) with
$\vartheta = \theta^1 + \theta^2$, and $\d \vartheta = 0$.

\section{Bicovariant differential calculus on the quantum group $GL_{p,q}(2)$}
\label{sec:GLpq(2)}
\setcounter{equation}{0}
The matrix quantum group $GL_{p,q}(2)$ is the algebra $\A$ generated by elements $a,b,c,d$
subject to the relations
\bea
    \begin{array}{l@{\quad}l@{\quad}l}
     a \, b = p \, b \, a & a \, c = q \, c \, a & b \, c = \frac{q}{p} \, c \, b \\
     b \, d = q \, d \, b & c \, d = p \, d \, c & a \, d = d \, a + (p - \frac{1}{q}) \, b \, c
    \end{array}
\eea
and supplied with a Hopf algebra structure \cite{Klim+Schm97}.
Its automorphisms are given by scalings
$(a,b,c,d) \mapsto (\alpha \, a , \beta \, b , \gamma \, c , \delta \, d )$ with
constants satisfying $\alpha \, \delta = \beta \, \gamma \neq 0$.
The relations defining a bicovariant differential calculus on a matrix
quantum group are usually expressed in terms of left-coinvariant Maurer-Cartan 1-forms:
\bea
 \left( \begin{array}{cc} \vartheta^1 & \vartheta^2 \\
   \vartheta^3 & \vartheta^4 \end{array} \right)
 := {\cal S} \left( \begin{array}{cc} a & b \\
   c & d \end{array} \right) \; \d \left( \begin{array}{cc} a & b \\
   c & d \end{array} \right)
\eea
using the antipode ${\cal S}$ (see \cite{MH92,MH+Reut93} for details).
Combined with (\ref{da-e_s}), this implies
\be
    \left( \begin{array}{c} \vartheta^1 \\ \vartheta^3 \end{array} \right)
  = {\cal S} \left( \begin{array}{cc} a & b \\ c & d \end{array} \right)
    \sum_{s=1}^4 \left( \begin{array}{c} e_s(a) \\ e_s(c) \end{array} \right) \theta^s
    \label{MC-e_s}
\ee
and a corresponding formula with $\vartheta^1, \vartheta^3, e_s(a), e_s(c)$
replaced by $\vartheta^2, \vartheta^4, e_s(b), e_s(d)$, respectively.
There are two distinguished bicovariant calculi on $GL_{p,q}(2)$ \cite{Malt90,MH92}.
For one of them, 1-forms satisfying (\ref{theta-phi}) have been derived in \cite{DMH04auto}.
The other is given by
\bea
 \begin{array}{l}
 \vartheta^1 \, a = {1 \over r} \, a \, \vartheta^1 \, ,  \quad
 \vartheta^1 \, b = b \, \vartheta^1  \, , \quad
 \vartheta^2 \, a = - (q - {1 \over p}) \, b \, \vartheta^1 + {1 \over p} \, a \, \vartheta^2
                     \nonumber \\
 \vartheta^2 \, b = {1 \over p} \, b \, \vartheta^2  \, , \quad
 \vartheta^3 \, a = {1 \over q} \, a \, \vartheta^3 \, , \quad
 \vartheta^3 \, b = - (p - {1 \over q}) \, a \, \vartheta^1
                    + {1 \over q} \, b \, \vartheta^3 \nonumber \\
 \vartheta^4 \, a = a \, [ {1 \over r} (1-r)^2 \, \vartheta^1 + \vartheta^4 ]
                      + {1-r \over r} \, b \, \vartheta^3  \, , \quad
 \vartheta^4 \, b = {1-r \over r} \, a \, \vartheta^2 + {1 \over r} \, b \, \vartheta^4
 \end{array}
\eea
where the remaining relations are obtained from the above by replacing $a$ by $c$ and $b$ by $d$
(see example 2 in section 6 of \cite{MH92}).
A lengthy calculation yields the 1-forms
\bea
  \begin{array}{l}
  \tilde{\theta}^1 = \vartheta^1 \, , \quad
  \tilde{\theta}^2 = - r \, d \, \vartheta^1 + c \, \vartheta^2 \, , \quad
  \tilde{\theta}^3 = a \, \vartheta^1 + b \, \vartheta^3 \, , \\
  \tilde{\theta}^4 = - p \, d \, a \, \vartheta^1 + c \, a \, \vartheta^2
                     - (p/q) \, b \, d \, \vartheta^3 + (p/q) \, b \, c \, \vartheta^4
                     \label{ttheta}
  \end{array}
\eea
which satisfy (\ref{theta-phi}) with the automorphisms determined by
\bea
   \begin{array}{ll}
   \tilde{\phi}_1(a,b,c,d) = (r^{-1} a , b , r^{-1} c , d ) \, , &
   \tilde{\phi}_2(a,b,c,d) = (r^{-1} a , q^{-1} b , p^{-1} c , d ) \\
   \tilde{\phi}_3(a,b,c,d) = (r^{-1} a , q^{-1} b , p^{-1} c , d ) \, , &
   \tilde{\phi}_4(a,b,c,d) = (r^{-1} a , q^{-2} b , q p^{-1} c , d ) \, .
   \end{array}
\eea
Multiplying the 1-forms $\tilde{\theta}^s$ by non-zero constants, arbitrary positive integer
powers of $b$, $c$, and arbitrary integer powers of the quantum determinant
${\cal D} := a \, d - p \, b \, c$, preserves the property (\ref{theta-phi}),
but with different automorphisms.
In this way we can achieve $\tilde{\phi}_2 \neq \tilde{\phi}_3$, for example.
(\ref{ttheta}) can only be solved for the Maurer-Cartan 1-forms if we pass over to
the algebra $\hat{\A}$ obtained from $\A$ by requiring $b$ and $c$ to have inverses.
Comparing the resulting expressions with (\ref{MC-e_s}) then shows that $e_s$, $s \in S$,
do not have values in $\A$, but rather in $\hat{\A}$. Let us write
\be
    \tilde{\theta}^s = b \, c \, \theta^s
\ee
with new (generalized) 1-forms $\theta^s$. The corresponding automorphisms are given by
\bea
   \begin{array}{ll}
   \phi_1(a,b,c,d) = (a , q p^{-1} b , q^{-2} c , r^{-1} d ) \, , &
   \phi_2(a,b,c,d) = (a , p^{-1} b , q^{-1} c , r^{-1} d ) \\
   \phi_3(a,b,c,d) = (a , p^{-1} b , q^{-1} c , r^{-1} d ) \, , &
   \phi_4(a,b,c,d) = (a , r^{-1} b , c , r^{-1} d )
   \end{array}
\eea
and (\ref{ttheta}) leads to
\bea
 \begin{array}{l}
 \vartheta^1 = b \, c \, \theta^1 \, , \quad
 \vartheta^2 = \frac{q}{p} \, b \, ( d \, \theta^1 + \theta^2 ) \, , \quad
 \vartheta^3 = - p \, a \, c \, \theta^1 + c \, \theta^3 \, , \\
 \vartheta^4 = - q \, a \, d \, \theta^1 - q \, a \, \theta^2
               + \frac{1}{p} \, d \, \theta^3 + \frac{q}{p} \, \theta^4
 \end{array}
\eea
The calculus is inner with (see \cite{MH92})
\be
    \vartheta = \frac{r}{1-r} ( r \, \vartheta^1 + \vartheta^4 )
  = \frac{q}{1-r} ( - r \, {\cal D} \, \theta^1 - r \, a \, \theta^2
    + d \, \theta^3 + q \, \theta^4 ) \; .
\ee
Comparison with (\ref{vartheta}) shows that $a_1 = -(qr/1-r) {\cal D}$,
$a_2 = -(qr/1-r) a$, $a_3 = (q/1-r) d$ and $a_4 = (q^2/1-r) {\bf 1}$ (where
$\bf 1$ is the unit in $\A$). This determines the twisted inner derivations
(\ref{twist-inner-deriv}). By rescalings of the $\theta^s$ we can
achieve somewhat simpler expressions.
Assuming again $b$ and $c$ to be invertible, the Maurer-Cartan 2-form relations
(see \cite{MH92}) translate to the considerably simpler relations
\bea
    \begin{array}{l}
    (\theta^1)^2 = (\theta^2)^2 = (\theta^3)^2 = (\theta^4)^2 = 0 \, , \\
    \theta^2 \, \theta^1 = - r \, \theta^1 \, \theta^2 \, , \quad
    \theta^3 \, \theta^1 = - \theta^1 \, \theta^3 \, , \quad
    \theta^3 \, \theta^2 = - \theta^2 \, \theta^3 \, , \\
    \theta^4 \, \theta^1 = - r \, \theta^1 \, \theta^4
        + (p - \frac{1}{q}) \, \theta^2 \, \theta^3 \, , \quad
    \theta^4 \, \theta^2 = - r \, \theta^2 \, \theta^4 \, , \quad
    \theta^4 \, \theta^3 = - r \, \theta^3 \, \theta^4 \, .
    \end{array}
\eea
The full differential calculus is inner, i.e. $\Xi = 0$.

It is quite surprising that for bicovariant calculi on $GL_{p,q}(2)$,
with the extension of the algebra by inverses of $b$ and $c$, the
bimodule of 1-forms splits into a direct sum of 1-dimensional bimodules,
so that we have the `automorphism structure' (\ref{theta-phi}).
We do not know yet whether this result extends to
bicovariant calculi on other matrix quantum groups (see \cite{Klim+Schm97},
for example).\footnote{Even for $GL_{p,q}(2)$ we do not have the full answer yet,
since there is a 1-parameter family of bicovariant calculi \cite{MH92,MH+Reut93}
and we were only able so far to treat special members.}
Is there a deeper relation with bicovariance? Let $\vartheta^s$ be the corresponding
\emph{left}-invariant Maurer-Cartan forms, so that
$\Delta_L \vartheta^s = 1 \otimes \vartheta^s$ with the left-coaction $\Delta_L$.
Writing $\theta^s = \sum_{s' \in S} \hat{a}^s{}_{s'} \; \vartheta^{s'}$, we find
\be
     \sum_{s' \in S} \Delta ( \hat{a}^s{}_{s'} ) \, 1 \otimes \vartheta^{s'} \, \Delta(a)
   = \Delta_L( \theta^s ) \, \Delta(a)
   = \Delta_L( \theta^s \, a )
\ee
where $\Delta$ is the coproduct. For all $a \in \A$, this has to be equal to
\be
    \Delta(\phi_s \, a) \, \Delta_L(\theta^s)
  = \Delta(\phi_s \, a) \, \sum_{s' \in S} \Delta( \hat{a}^s{}_{s'}) \, 1 \otimes \vartheta^{s'}
    \; .
\ee
Additional conditions arise from \emph{right}-covariance of the calculus.
Together with the knowledge of the corresponding automorphism groups, these conditions
restrict the elements $\hat{a}^s{}_{s'}$ and should lead to the possible choices of
1-forms $\theta^s$ or corresponding obstructions. This program has still to be carried out.

\section{Differential geometry in terms of the basis $\theta^s$}
\label{sec:diffgeo}
\setcounter{equation}{0}
An automorphism $\phi$ of $\A$ is called \emph{differentiable} with respect to a
differential calculus $(\hat{\Omega}(\A),\d)$ if it extends to a homomorphism of
$\hat{\Omega}(\A)$ such that $\phi \circ \d = \d \circ \phi$.
In the above examples, the $\phi_s$ are indeed differentiable.
For the $q$-plane WZ-calculus, for example, one finds $\phi_s \, \theta^{s'} = \theta^{s'}$.
Let us now assume that $\{ \phi_s , \phi_s^{-1} \, | \, s \in S \}$ are differentiable.
Then we can introduce a \emph{left $\A$-linear} tensor product basis:
\be
      \theta^s \otimes_L \theta^{s'}
   := \theta^s \otimes_{\A} \phi_s^{-1} \theta^{s'} \quad
  \Longrightarrow \quad a \, \theta^s \otimes_L ( a' \, \theta^{s'} )
   = a \, a' \, \theta^s \otimes_L \theta^{s'}  \; .
\ee
In particular, a \emph{metric} should then be an element\footnote{See \cite{DMH03grpl2}
for arguments why a metric should be defined in this way rather than as an element of
$\hat{\Omega}^1 \otimes_{\A} \hat{\Omega}^1$.}
\be
   \mathbf{g} = \sum_{s,s' \in S} g_{s,s'} \; \theta^s \otimes_L \theta^{s'}
\ee
with $g_{s,s'} \in \A$.
Let $\nabla$ be a \emph{linear connection}, i.e. a connection on
$\hat{\Omega}^1$. Then
\be
    \nabla \omega = \vartheta \otimes_{\A} \omega
    - \sum_{s' \in S} \theta^{s'} \otimes_{\A} \V_{s'} \, \omega
    \qquad \forall \, \omega \in \hat{\Omega}^1
\ee
with `parallel transport operators' $\V_s$ satisfying
\be
    \V_s ( a \, \omega ) = (\phi_s^{-1} \, a) \; \V_s \, \omega  \; .
\ee
Furthermore,
\be
    \V_s ( \omega \otimes_L \omega' ) := (\V_s \, \omega) \otimes_L (\V_s \, \omega' )
    \qquad \forall \omega, \omega' \in \hat{\Omega}^1
\ee
determines a connection $\nabla$ on
$\hat{\Omega}^1 \otimes_L \hat{\Omega}^1$.
So we can define \emph{metric compatibility}:
\be
     \nabla \mathbf{g} = 0  \quad \Longleftrightarrow \quad
     \V_s(\mathbf{g}) = \mathbf{g} \quad \forall s \in S \; .
\ee
A torsion-free metric-compatible linear connection should be called a {\em Levi-Civita
connection}. For differential calculi associated with (a subclass of) Cayley graphs
(`group lattices') of finite groups, the corresponding generalized (pseudo-)
Riemannian geometry has been developed in \cite{DMH03grpl1,DMH03grpl2}.
This includes (pseudo-) Riemannian geometry of hyper-cubic lattices as a special case.
The above scheme can be applied as well to introduce analogs of (pseudo-) Riemannian
structures on quantum planes and quantum groups.

\end{document}